\documentclass[graybox]{svmult}
\usepackage{mathptmx}       
\usepackage{makeidx}         
\usepackage{graphicx}        
\usepackage{multicol}        
\usepackage[bottom]{footmisc}
\usepackage{url}

\usepackage{algorithm}
\usepackage{algpseudocode}
\usepackage{amsfonts}
\usepackage{amsmath}
\usepackage{amssymb}
\usepackage{mathtools}
\usepackage{wrapfig}
\usepackage{csquotes}
\usepackage{floatflt}
\usepackage[font=small,skip=0pt]{caption}
\usepackage{adjustbox}

\graphicspath{{./Graphics}}

\begin{document}

\title{A global optimum-informed greedy algorithm for A-optimal experimental design }
\titlerunning{Global optimum-informed greedy OED}
\author{Christian Aarset} 
\institute{University of G\"ottingen - Institute for Numerical and Applied Mathematics, Lotzestr. 16-18, D-37083 G\"ottingen, Germany  
} 

\maketitle

\abstract{Optimal experimental design (OED) concerns itself with identifying ideal methods of data collection, e.g.~via sensor placement. The \emph{greedy algorithm}, that is, placing one sensor at a time, in an iteratively optimal manner, stands as an extremely robust and easily executed algorithm for this purpose. However, it is a priori unclear whether this algorithm leads to sub-optimal regimes. Taking advantage of the author's recent work on non-smooth convex optimality criteria for OED, we here present a framework for rejection of sub-optimal greedy indices, and study the numerical benefits this offers.}


\keywords{optimal experimental design, inverse source problem, non-smooth convex optimality criteria, greedy algorithm, A-optimality, finite elements}
\\
{{\bf MSC2020:} 62K05, 62F15, 35R30, 65K10.}


\def\R{\mathbb{R}}
\def\C{\mathbb{C}}
\def\N{\mathbb{N}}
\def\F{\mathbb{F}}
\def\M{\mathbb{M}}

\def\mc{\mathbf{m}}

\def\Cc{\mathcal{C}}
\def\Jc{J}
\def\Fc{\mathcal{F}}
\def\Nc{\mathcal{N}}

\def\Om{\Omega}
\def\Ls{L^2(\Omega_0)}
\def\Ld{L^2(\Omega_1,\C)}
\renewcommand{\Re}{\mathop{\mathrm{Re}}\nolimits}
\renewcommand{\Im}{\mathop{\mathrm{Im}}\nolimits}
\renewcommand{\d}{\,\mathrm{d}}

\def\tr{\mathop{tr}\nolimits}
\newcommand{\ws}{w^*}
\newcommand{\wb}{w^\#}
\def\prior{\text{prior}}
\def\post{\text{post}}
\def\mp{\mc_{\prior}}
\def\mpst{\mc_{\post}}
\def\Cp{\Cc_{\prior}}
\def\Cpi{\Cc^{-1}_{\prior}}
\def\Cpst{\Cc_{\post}}

\newcommand{\blue}[1]{\textcolor{blue}{#1}}
\section{Introduction} 

Optimal experimental design (OED) can be seen as the field of identifying \emph{designs} $w$ allowing for the best reconstruction of unknown parameters $x\in X$ in some ambient space $X$, given that $x$ can only be measured indirectly by some $w$-dependent forward map $F_w$, i.e.~one only has access the noisy, design-dependent data $g_w$ given by
\begin{equation}\label{eq:inv}
    g_w = F_w f + \epsilon, \quad \epsilon\sim\Nc(0,\sigma^2I), \quad \sigma>0
\end{equation}
in the case of Gaussian white noise. A typical example of the effect of the design on the experiment is the situation where $F_w = M_w F$, where $F:X\to\R^m$ and $m\in\N$ are fixed -- ubiquitously, $F$ is a composition of a finite observation operator and a partial differential equation (PDE) solution operator \cite{Alexanderian} -- and $M_w\in\R^{m\times m}$ is the diagonal matrix with $w\in\R^m$ on the diagonal. Imposing $w\in\{0,1\}^m$, the design $w$ acts as a mask on the data, representing \emph{sensor placement} -- the quantity $m$ being the number of candidate locations where the experimenter might elect to place a sensor\footnote{\cite{Aarset} details the situation where each entry of $F_wf$ is a vector of measurements, e.g.~each sensor observing at multiple time points or observing multiple frequencies.}, while each index $k\in\N$, $k\leq m$ corresponds to the experimenter's choice of either making a measurement of $Ff$ in the $k$-th candidate location, in which case $w_k:=1$, or not to make it, in which case $w_k:=0$.

If, due to e.g.~budget or power constraints, one can place only $m_0<m$ sensors, then this leads to the \emph{sensor placement problem} -- identifying the best selection of $m_0$ out of $m$ candidate locations to use. In order to determine one design $w$ as better than another, one fixes as objective a design criterion $\Jc:\R^m\to\R$, mapping designs $w$ to a measure $\Jc(w)$ of the quality of the reconstruction of $f$. A number of such design criteria exist, including but not limited to A-optimality, D-optimality and expected information gain; \cite{Alexanderian, Pukelsheim} offer comprehensive overviews. In this article, we will solely concern ourselves with A-optimality, as detailed in the next section.


Various methods have been proposed towards solving the sensor placement problem for infinite-dimensional inverse problems governed by PDEs, also by the author \cite{Aarset}. In the present article, we will focus on arguably the most elementary algorithm; namely, the  \emph{greedy} algorithm\footnote{Various interesting adaptations of the greedy algorithm have been proposed, e.g.~\cite{WuChenGhattas}.}. Its advantage is its simplicity, requiring only approximately $mm_0-m_0^2/2$ evaluations of the objective $\Jc$. However, it can easily enter sub-optimal regimes. The main contribution of this article will be demonstrating how the optimality criteria developed in \cite{Aarset} can be used to adaptively correct the output of the greedy algorithm. 


%

\section{A-optimal designs}

We here briefly introduce the A-optimal criterion; for further details, we refer to \cite{Ucinski}. In the Bayesian setting, the linear inverse problem \eqref{eq:inv} given prior distribution $f\sim\Nc(\mp,\Cp)$ has explicit, design-dependent posterior distribution
\begin{align}\label{eq:bayesian}
    f \mid g_w & \sim \Nc(\mpst,\Cpst), \nonumber \\
    \mpst(w) & := \mp + \Cpst F_w^*\Sigma^{-1}\left(
        g_w - F_w \mp
    \right), \\
    \Cpst(w) & := \left(
        F_w^*\Sigma^{-1}F + \Cpi
    \right)^{-1}, \nonumber
\end{align}
see \cite{Stuart}. The A-optimal objective $\Jc:\R^m\to\R$ is given as $\Jc(w):=\tr(\Cpst(w))$, which by Mercer's theorem can be seen as proportional to the pointwise variance in the reconstruction, which an A-optimal design $w^*$ thus minimises. 

\section{Optimality for sensor placement}

Given the above, one can cast the problem of finding the A-optimal design $\wb$ using exactly $m_0$ sensors as the constrained optimisation problem of determining
\begin{equation}\label{eq:OED0}
\wb\in\mathop{\mathrm{argmin}}_{w\in K_0}\,\Jc(w), \qquad K_0:=\{
    w\in\{0,1\}^m \,\mid\, \|w\|_0\leq m_0
\}
\end{equation}
with the zero-\enquote{norm} $\|w\|_0:=\#\{k\in\N \,\mid \,k\leq m, \,w_k\neq 0\}$. Since this is a binary, non-convex optimisation problem, it is difficult to treat exactly. As in \cite{Aarset}, we will instead consider the so-called $1$-relaxed problem of finding the \emph{global optimal design}
\begin{equation}\label{eq:OED}
\ws\in\mathop{\mathrm{argmin}}_{w\in K_1}\,\Jc(w), \qquad K_1:=\{
    w\in\R^m \mid 0\leq w\leq 1, \, \sum_{k=1}^mw_k\leq m_0
\},
\end{equation}
which, as argued in \cite{Aarset}, is a convex optimisation problem and thus is generally solvable. While $\ws$ will typically be non-binary, with multiple indices taking values between $0$ and $1$, it serves as a lower bound on the A-optimality of the binary optimal design $\wb$, in the sense that $\Jc(\ws)\leq\Jc(\wb)$, since necessarily $\wb\in K_0\subseteq K_1$. Moreover, a key contribution of \cite{Aarset} was employing first-order optimality criteria for convex constrained problems to show that the non-binary global optimum $\ws$ \emph{does}, in fact, contain a surprising number of indices exactly equal to $0$ or $1$.

\begin{theorem}[{\cite[Thm.~2]{Aarset}}]Given $m_0\leq m$ and $w\in K_1$, assume (reordering if necessary) that the indices $k$ of $w$ are ordered so that
\[
    \nabla J(w)_1 \leq \nabla J(w)_2 \leq \ldots \leq \nabla J(w)_m.
\]

Then $w=\ws$ if and only $w_k=1$ for all $k$ satisfying $\nabla J(w)_k < \nabla J(w)_{m_0+1}$, $w_k=0$ for all $k$ satisfying $\nabla J(w)_k > \nabla J(w)_{m_0}$ and $\sum_{k=1}^mw_k=m_0$.
\label{thm:optimality}
\end{theorem}

Theorem \ref{thm:optimality} is remarkable, in that it provides a convex constrained optimality criterion that can be verified explicitly via access to the gradient of the objective functional. Moreover, \cite{Aarset} demonstrates that the generic situation for $\ws$ is that a small number (possibly zero) of \emph{dominant} indices $k$ satisfy $\ws_k=1$ exactly, a large number of \emph{redundant} indices $k$ satisfy $\ws_k=0$ exactly, and remaining indices $k$ satisfy $\Jc(\ws)_k=\Jc(\ws)_{m_0}$ exactly, with no statement on the value of $w_k\in[0,1]$. A strength of this formulation is that it is not sensitive to numerical error, as the ordering of the gradient can be used to clearly distinguish whether $\ws_k=0$ or $\ws_k=1$, as opposed to only numerically approximating these values. This means $\ws$ can be numerically found by standard algorithms, such as \texttt{scipy.optimize.minimize} \cite{scipy.optimize.minimize}, with minimal concern for numerical error. \cite{Aarset} moreover lays out how the evaluation, gradient and Hessian of $\Jc$ can be computed extremely cheaply for the A-optimal objective, requiring no PDE solves or trace estimation.


\section{Greedy algorithms}

Given a number $m$ of candidate locations, the greedy algorithm approximates the optimal binary design $\wb$ for \emph{each} $m_0<m$ iteratively, by first testing the objective value $\Jc(w)$ for every configuration using exactly one sensor, then repeating, each time adding the single sensor that improves the previously found design the most. Explicitly, this leads to Algorithm \ref{alg:greedy}, where $e^k\in\R^m$ denotes the $k$-th unit vector, $e^k_l=\delta_{k=l}$. For each $m_0$, one then approximates $\wb$ by the output $w^{m_0}$. 

\vspace{-0.45cm}
\begin{algorithm}
\caption{Greedy OED}\label{alg:greedy}
\begin{algorithmic}
\Require $m_0:=0$, $w^0:=0$
\While{$m_0<m$}
\State Compute $k'\in\mathop{\mathrm{argmin}}_{k\in I_{m_0}}\Jc(w^{m_0}+e^k)$, where $I_{m_0}:=\left\{k\in\N,\,k\leq m \mid w^{m_0}_k=0\right\}$
\State $w^{m_0+1} \gets w^{m_0}+e^{k'}$ and $m_0 \gets m_0+1$
\EndWhile
\State \Return Each approximate binary optimal design $w^{m_0}$
\end{algorithmic}
\end{algorithm}
\vspace{-0.5cm}

Algorithm \ref{alg:greedy} finds each binary approximate optimal design $w^{m_0}\in\{0,1\}^m$ in a \emph{nested} manner. While this dramatically reduces search complexity compared to a naive full binary search, going from at most $\binom{m}{m_0}$ to scaling at worst quadratically in $m_0$, it is not a priori clear whether each $w^{m_0}$ is a good approximation of the optimal binary design $\wb$, nor is it clear \emph{when} the approximation deteriorates, whether it later improves, or if one permanently enters a sub-optimal regime. 

The optimality criterion Theorem \ref{thm:optimality} contributes two significant improvements to the above situation. Firstly, it allows for comparison between the A-optimal objective $\Jc(w^{m_0})$ and its lower bound $\Jc(\ws)$, giving some indication of the relative quality of $w^{m_0}$. More pressingly, the presence of \emph{redundant} indices $k$, i.e.~indices such that $\ws_k=0$ exactly, suggest that the relative contribution of these sensors is negligible. If $w^{m_0}_k=1$ for a redundant index $k$, one may thus be tempted to assume that this sensor is contributing to getting \enquote{stuck} in a sub-optimal regime. One may then attempt to build a global optimum-informed greedy design $w_*^{m_0}$ by setting redundant indices to $0$, then re-applying the greedy algorithm to replace these lost sensors. Put together, this leads to the modified greedy algorithm:

\vspace{-0.35cm}
\begin{algorithm}
\caption{Global optimum-informed greedy OED}\label{alg:fgreedy}
\begin{algorithmic}
\Require $m_0>0$ small
\State Compute $w^{m_0}$ via Algorithm \ref{alg:greedy}, $w^{m_0}_*\gets w^{m_0}$
\While{$m_0<m$}
\State Compute $\ws$ solving \ref{eq:OED}
\For{$1\leq k\leq m$}
\If{$\ws_k=0$}
\State $(w_*^{m_0})_k\gets 0$
\EndIf
\EndFor
\State Update $w_*^{m_0}$ by Algorithm \ref{alg:greedy} until it has $m_0$ non-zero entries, $m_0\gets m_0+1$
\EndWhile
\State \Return Each global optimum-informed approximate binary optimal design $w_*^{m_0}$
\end{algorithmic}
\end{algorithm}
\vspace{-0.5cm}

\section{Numerical experiments}


We demonstrate our modified greedy algorithm to a numerical example. To allow for comparison, we re-use the setting of \cite{Aarset}. Explicitly, define $\Omega_0:=B_{0.35}(0)$, $\Omega_1:=B_1(0)\setminus\bigcup_{i=1}^3S_i$, with three rectangular sound-hard scatterers as depicted in the upcoming Figures. With $m=334$, the sensor locations $(x_k)_{k=1}^m\in\Omega_1$ represent the intersection of a uniform circular grid with $\Omega_1$.

Given a source $f\in L^2(\Omega_0)$, let $u_{\omega}\in H^2(\Omega_1)$ solve the Helmholtz equation with impedance boundary on $\partial B_1(0)$ and wave number $\omega$, see \cite{Aarset, ColtonKress}. In each sensor, real and complex pointwise measurements of $u_\omega$ are made for seven frequencies $\{\omega_1,\ldots,\omega_7\}=\{20, 25, 30, 35, 40, 45, 50\}$. Thus, the parameter-to-observable map $F_w:L^2(\Omega_0)\to\R^{14m}$ is for each $k\in\N$, $k\leq m$ given as $(Ff)_{k+(0:13)m}:=w_k(\Re u_{\omega_1}(x_k),\Im u_{\omega_1}(x_k),\ldots,\Re u_{\omega_7}(x_k),\Im u_{\omega_7}(x_k))\in\R^{14}$. Discretisation is carried out via the NGSolve package \cite{Schoberl}. The Helmholtz solutions $u_\omega$ on $\Omega_1$ are discretised via a complex second-order FEM space employing $n=14641$ degrees of freedom; for detailed treatment on FEM discretisation for the Bayesian setting, see \cite{BuithanhGhattasMartinStadler}. Similarly, the prior covariance $\Cp$ was densely defined on $L^2(\Omega_0)$ as two-times application of the solution operator mapping the source $f$ to the solution $u$ of the Laplacian with Robin boundary condition, then discretised in the FEM. Measurement noise level $\sigma^2$ in \eqref{eq:inv} and \eqref{eq:bayesian} was chosen proportionally to $1\%$ of the average pointwise variance of $F(f_s)$ over $10^4$ samples $f_s$ drawn from the prior distribution $\Nc(0,\Cp)$; see \cite{BuithanhGhattasMartinStadler}.

\newlength{\figsize}
\setlength{\figsize}{0.265\textwidth}

\newlength{\negsize}
\setlength{\negsize}{3pt}

\begin{figure}
\vspace{-0.5cm}
\centering
\includegraphics[trim={19.8cm 5cm 19.8cm 6cm},clip,width=\figsize,keepaspectratio]{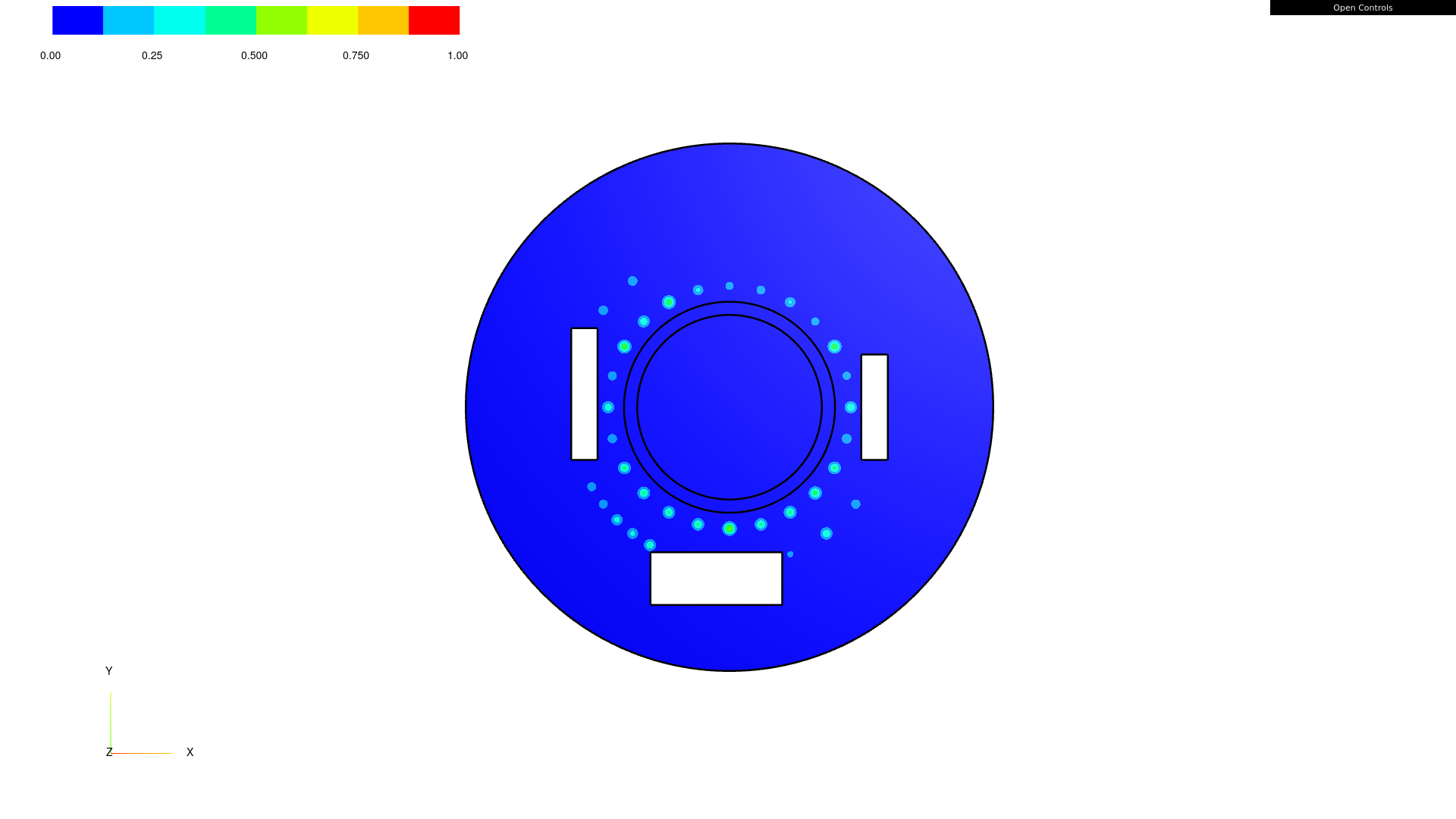}\hspace{-\negsize}
\includegraphics[trim={19.8cm 5cm 19.8cm 6cm},clip,width=\figsize,keepaspectratio]{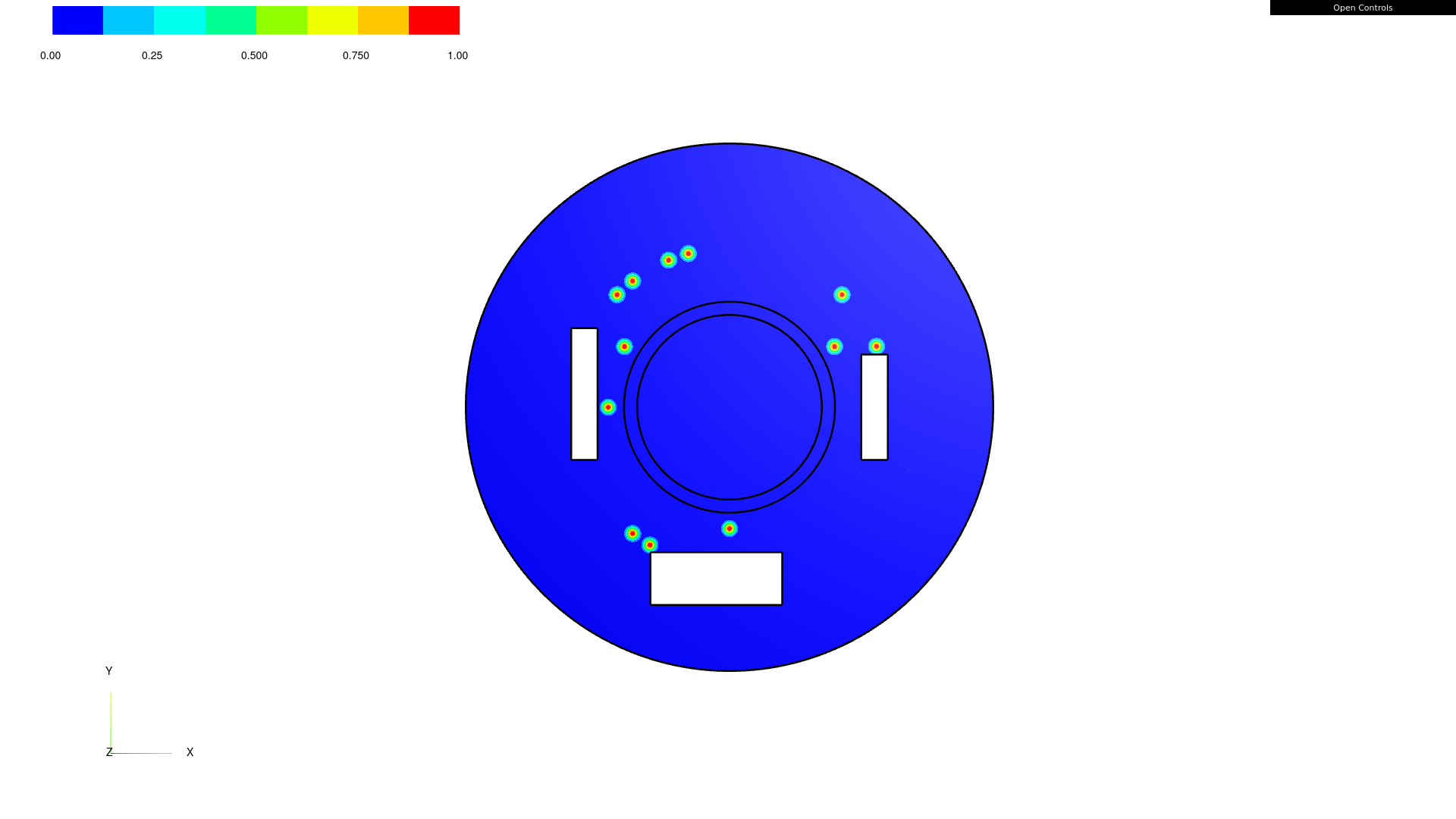} \hspace{-\negsize}
\includegraphics[trim={19.8cm 5cm 19.8cm 6cm},clip,width=\figsize,keepaspectratio]{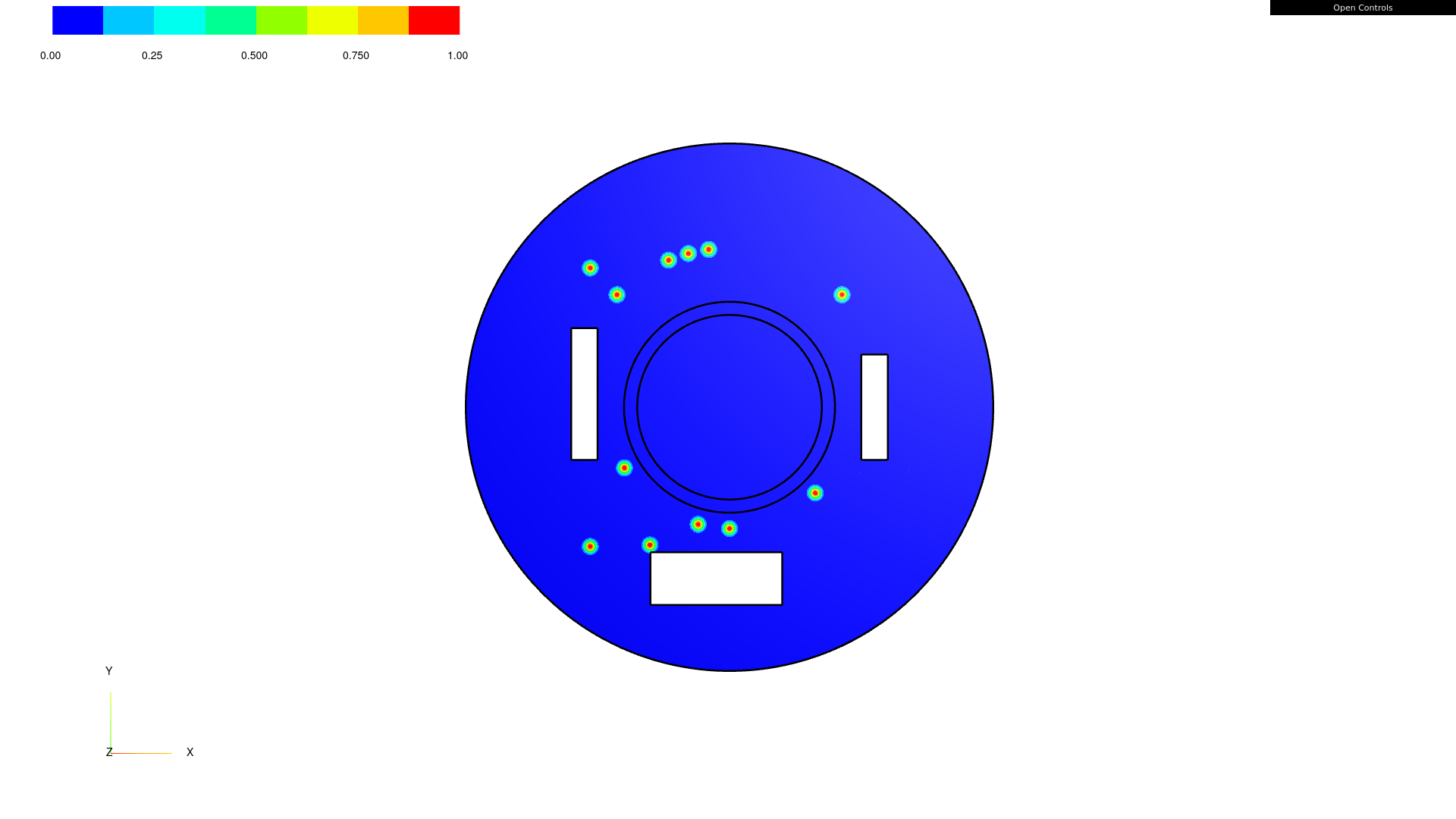} \\[-0.5ex]
\includegraphics[trim={19.8cm 5cm 19.8cm 6cm},clip,width=\figsize,keepaspectratio]{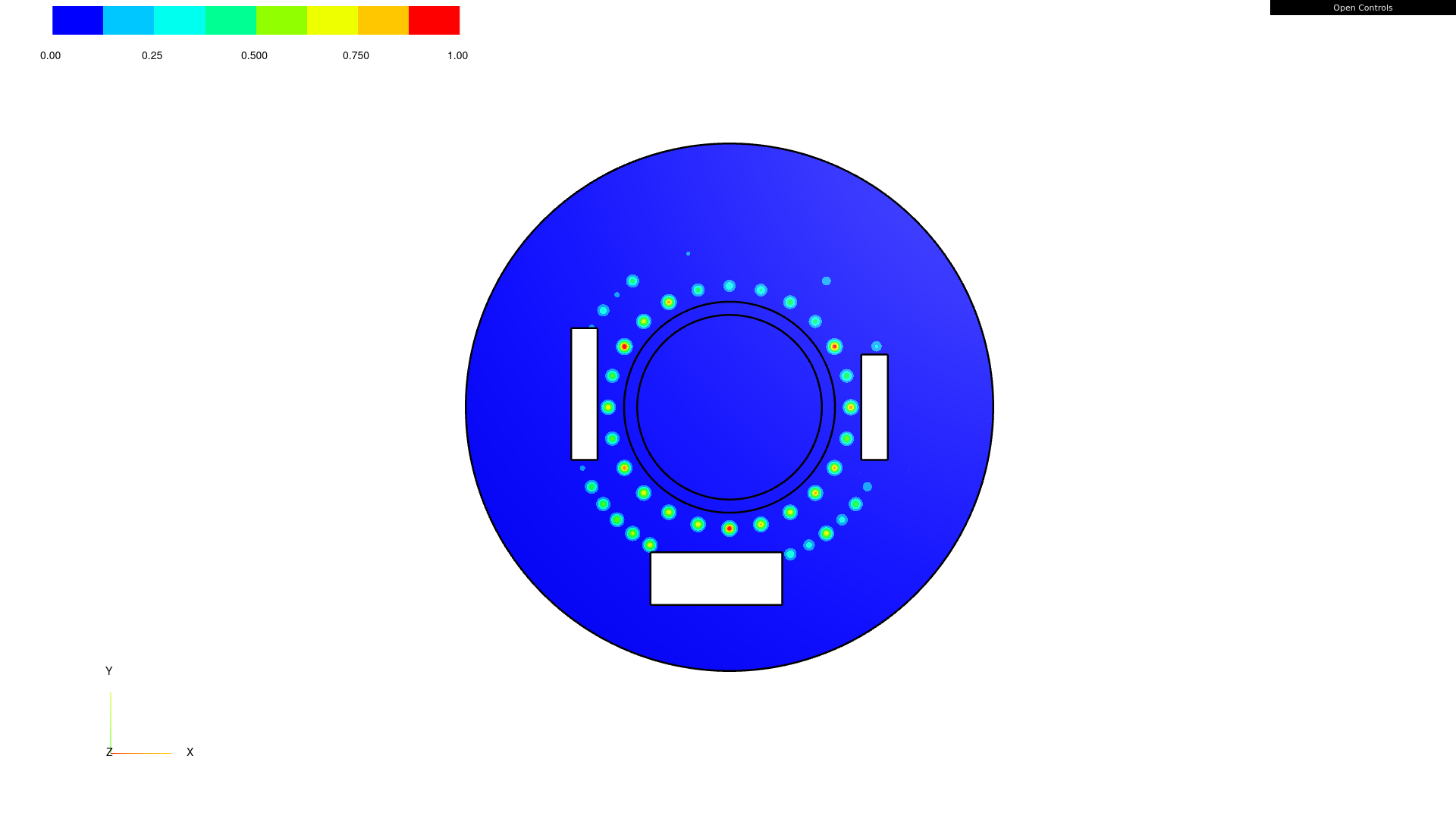}\hspace{-\negsize}
\includegraphics[trim={19.8cm 5cm 19.8cm 6cm},clip,width=\figsize,keepaspectratio]{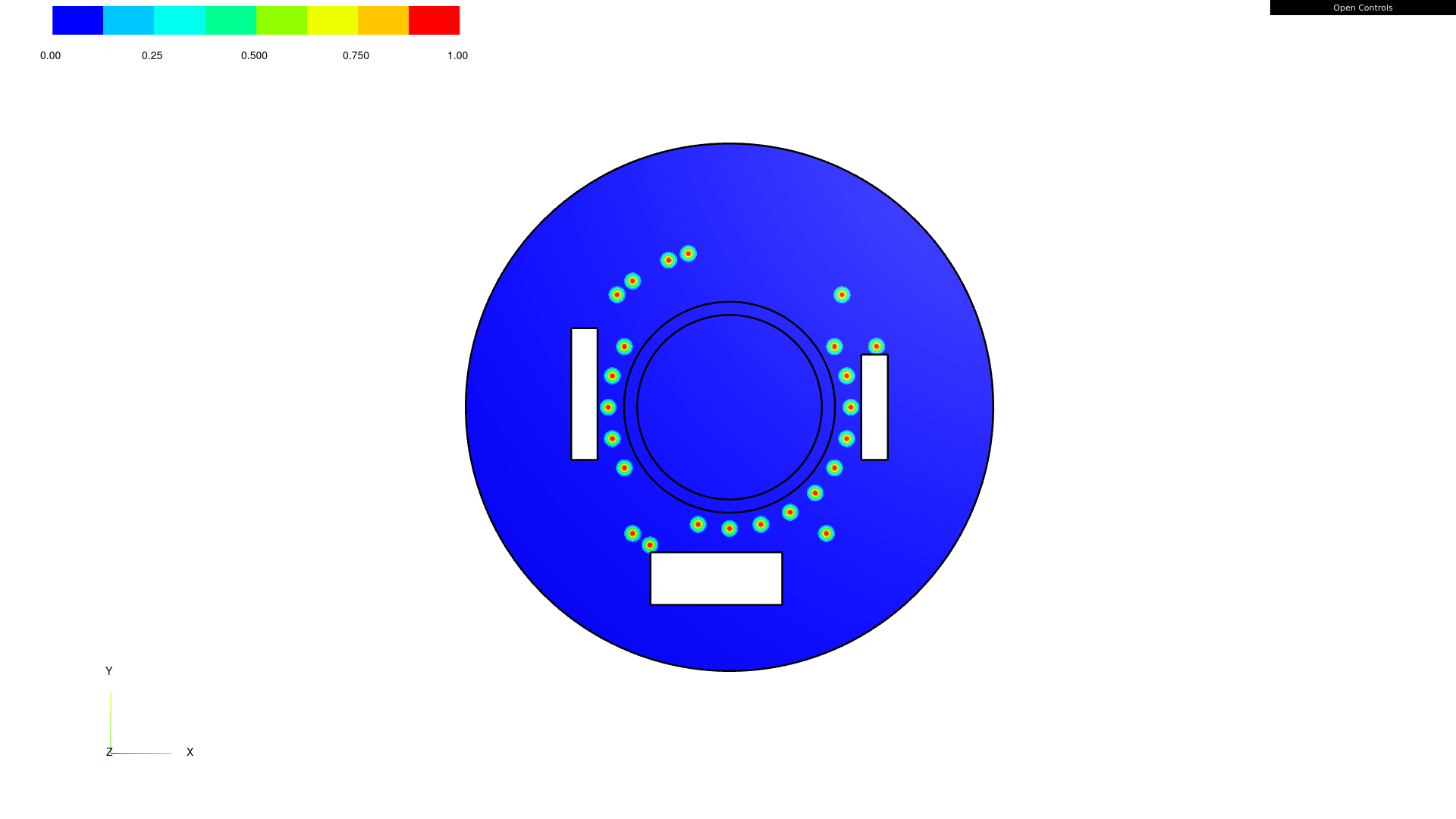} \hspace{-\negsize}
\includegraphics[trim={19.8cm 5cm 19.8cm 6cm},clip,width=\figsize,keepaspectratio]{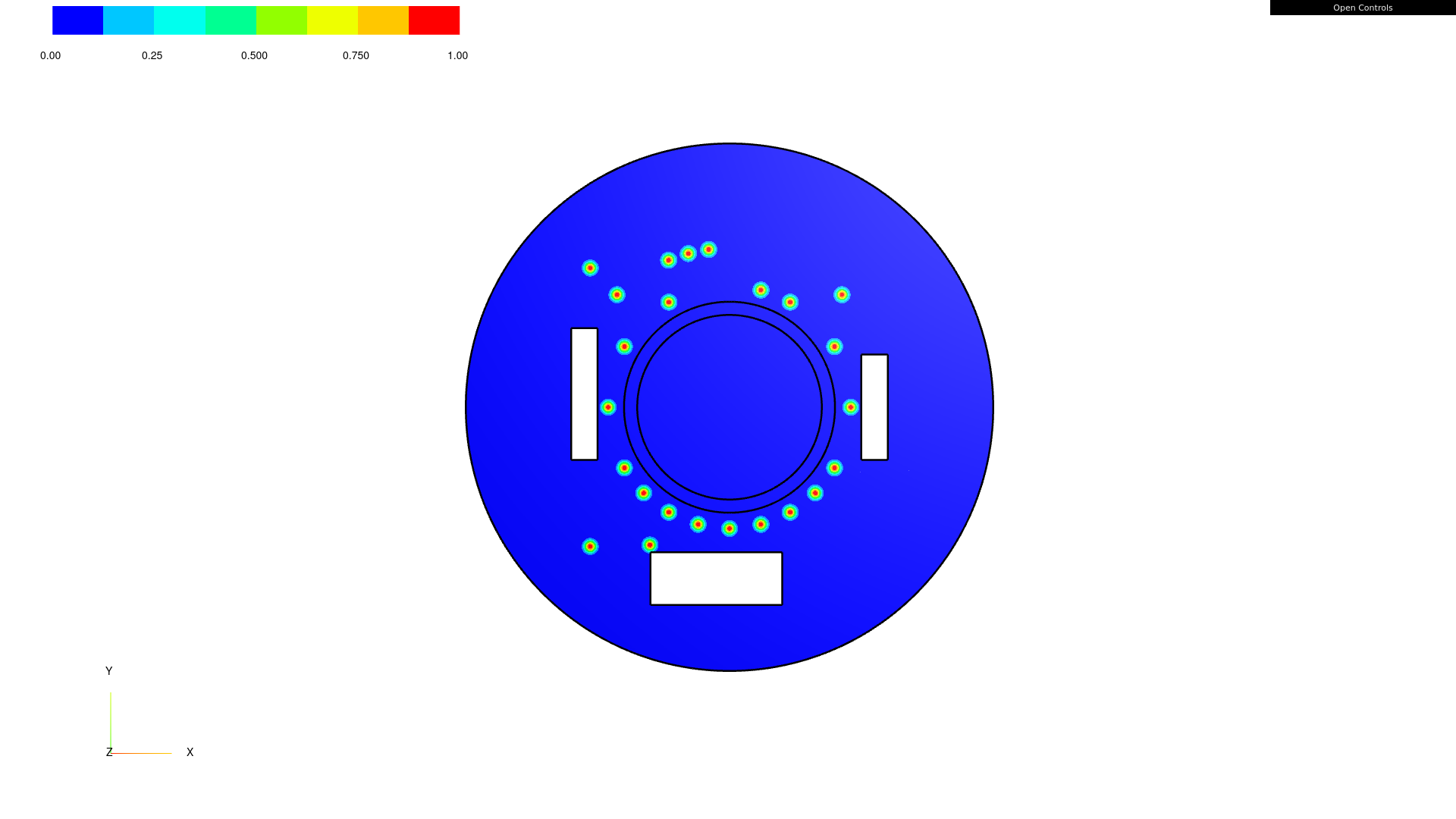} \\[-0.5ex]
\includegraphics[trim={19.8cm 5cm 19.8cm 6cm},clip,width=\figsize,keepaspectratio]{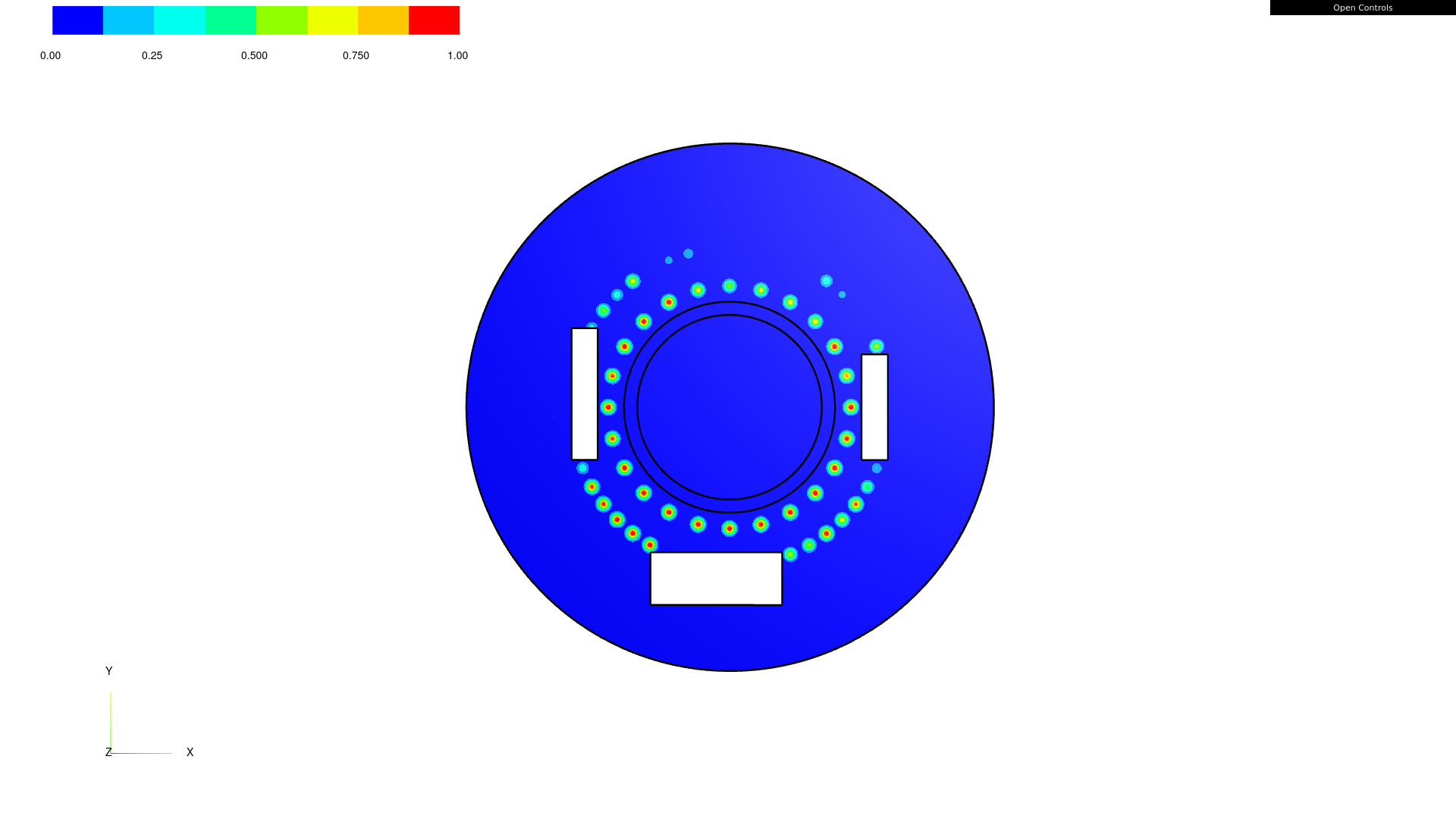}\hspace{-\negsize}
\includegraphics[trim={19.8cm 5cm 19.8cm 6cm},clip,width=\figsize,keepaspectratio]{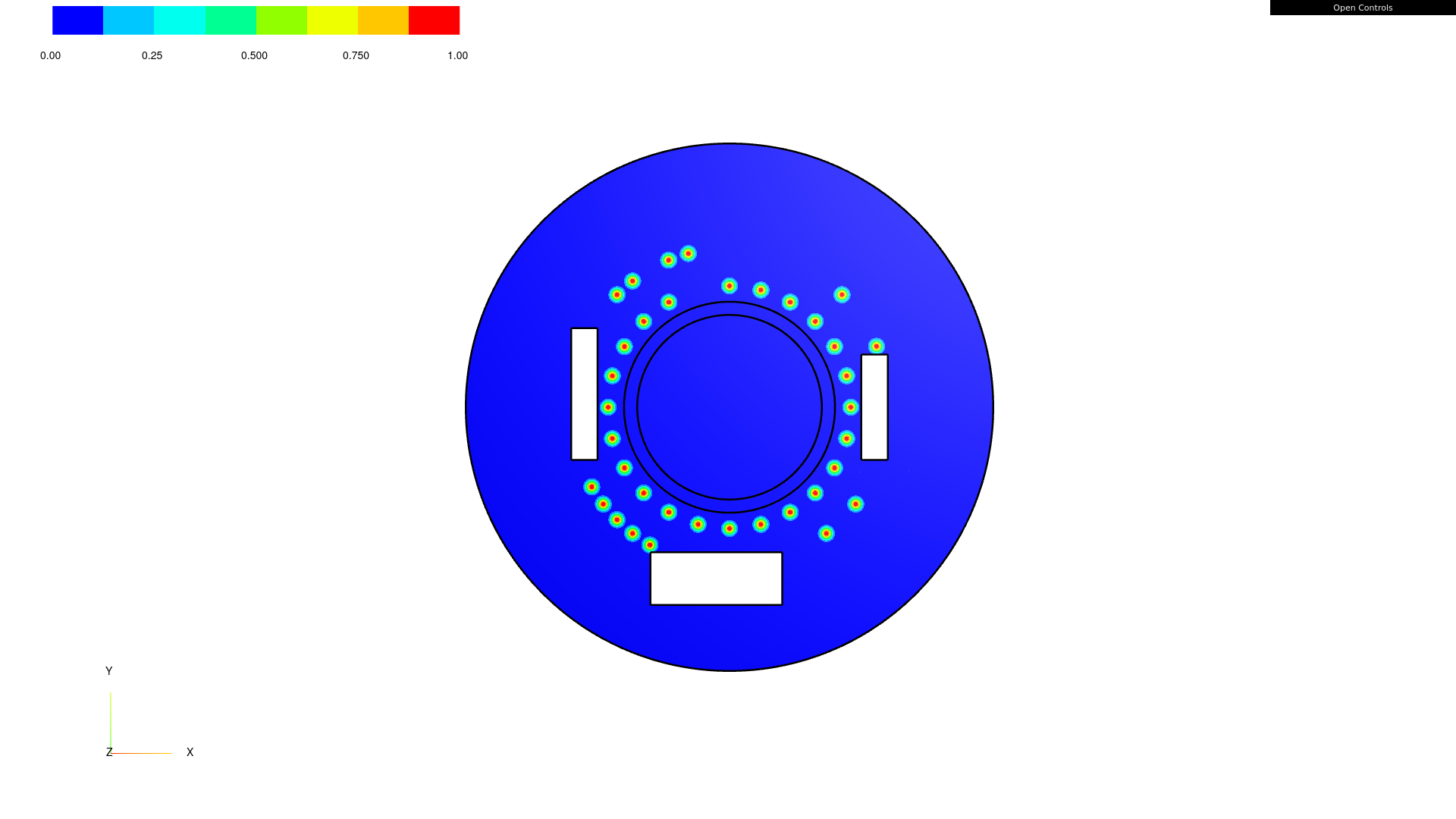} \hspace{-\negsize}
\includegraphics[trim={19.8cm 5cm 19.8cm 6cm},clip,width=\figsize,keepaspectratio]{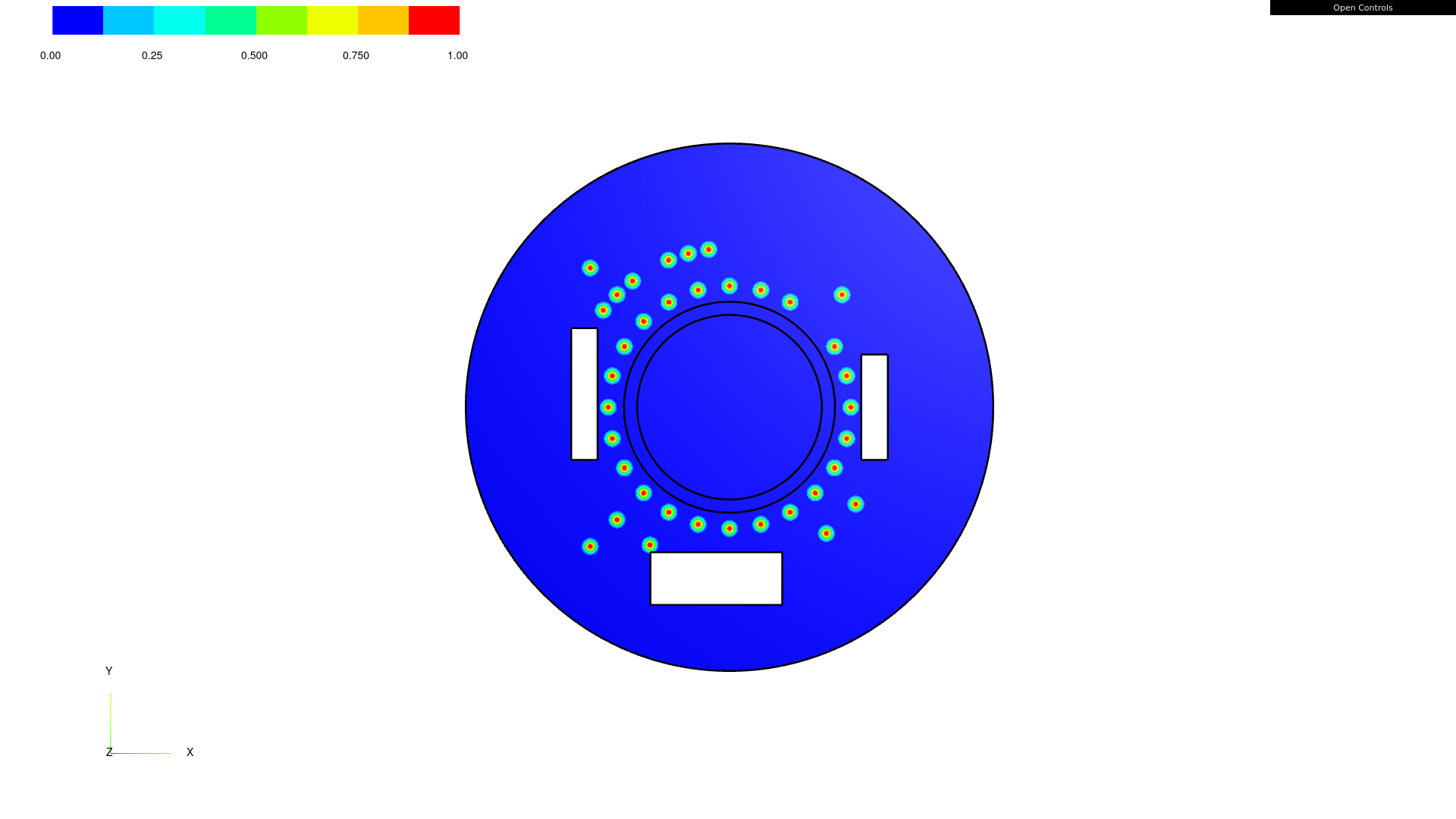}
\caption{Left to right: Globally optimal designs $\ws$, global optimum-informed greedy designs $w^{m_0}_*$, greedy designs $w^{m_0}$. Top to bottom: Designs for $m_0=12$, $24$, $36$.}
\vspace{-0.3cm}
\label{fig:first}
\end{figure}

\begin{figure}
\vspace{-1cm}
\centering
\includegraphics[width=0.65\textwidth,keepaspectratio]{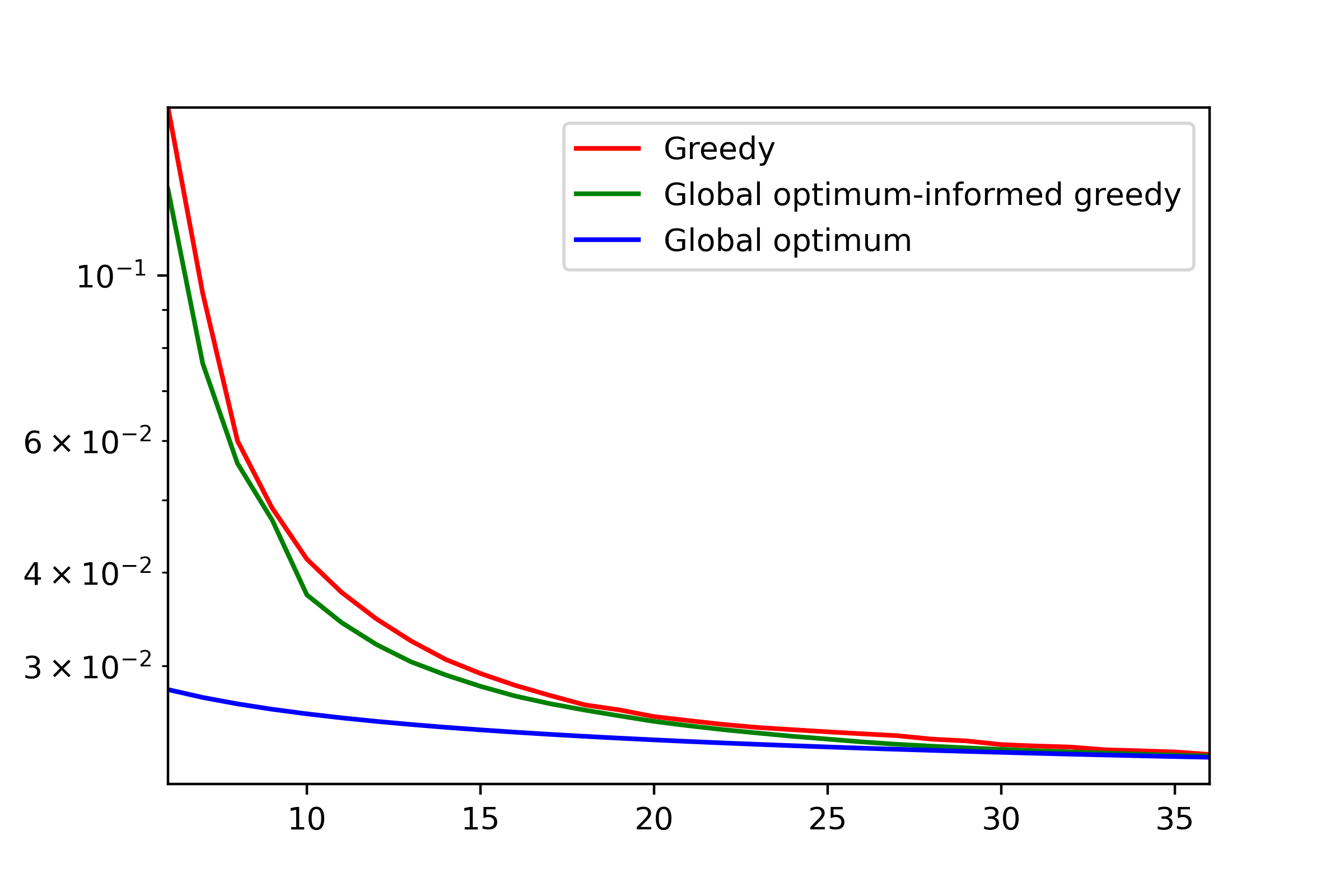}
\caption{A-optimal objectives $\Jc$ for the standard greedy sequence (red), global optimum-informed greedy sequence (green) and non-binary global optimum sequence (blue).}
\label{fig:objectives}
\vspace{-1cm}
\end{figure}

\begin{figure}
\begin{minipage}[t]{\linewidth}\vspace{0pt}
\setlength{\figsize}{0.265\textwidth}
\centering
\includegraphics[trim={19.8cm 5cm 19.8cm 6cm},clip,width=\figsize,keepaspectratio]{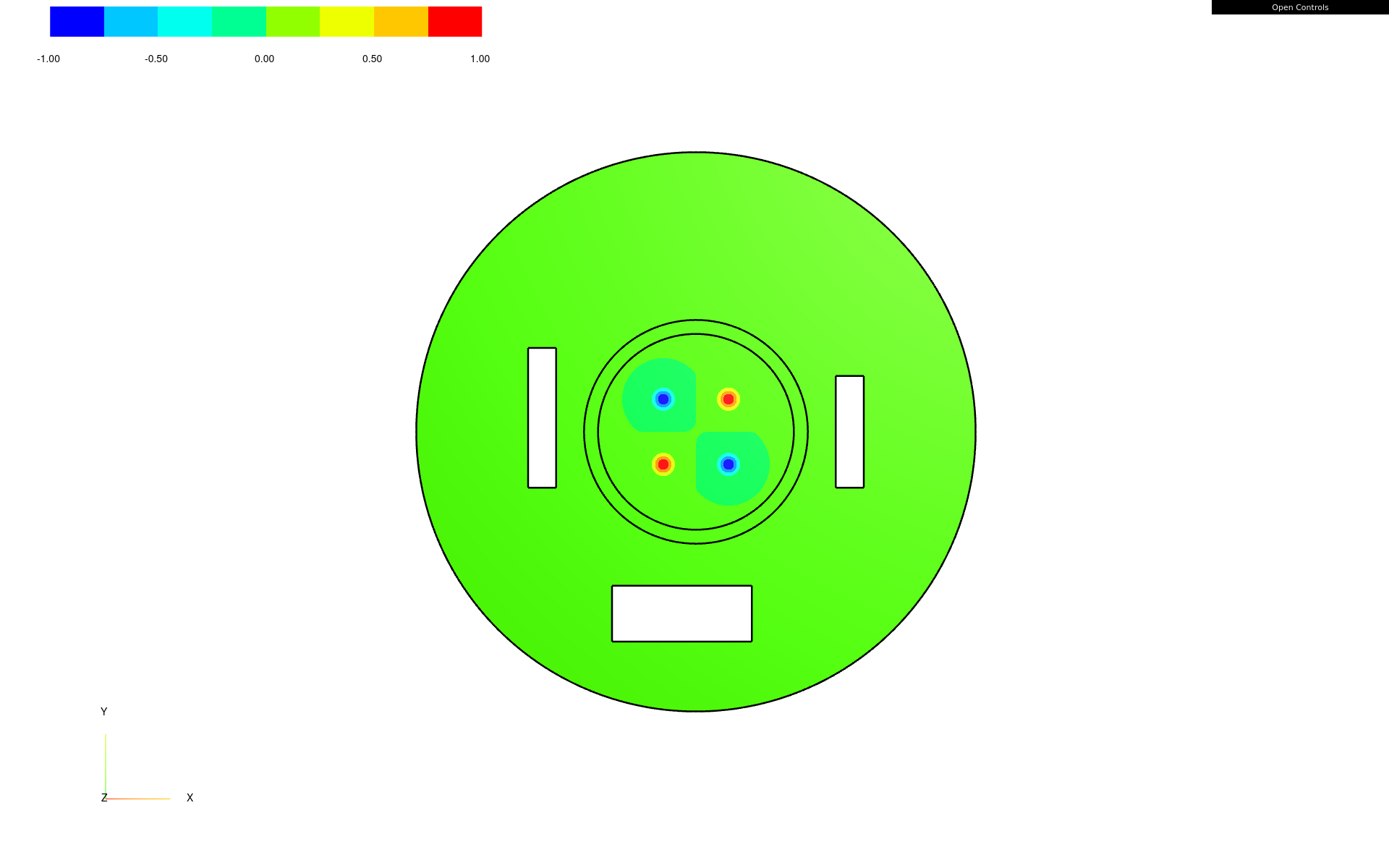}
\includegraphics[trim={19.8cm 5cm 19.8cm 6cm},clip,width=\figsize,keepaspectratio]{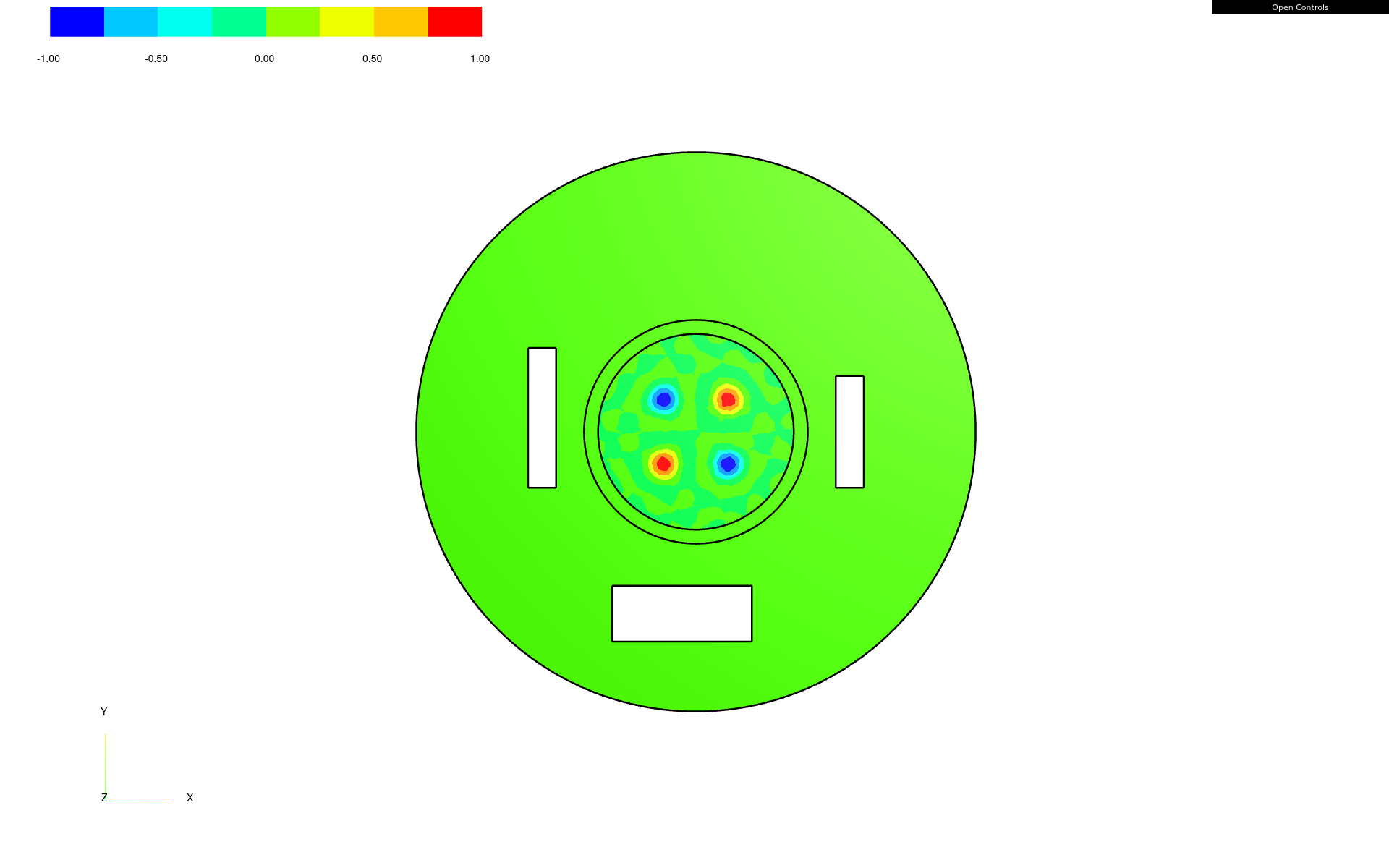}
\includegraphics[trim={19.8cm 5cm 19.8cm 6cm},clip,width=\figsize,keepaspectratio]{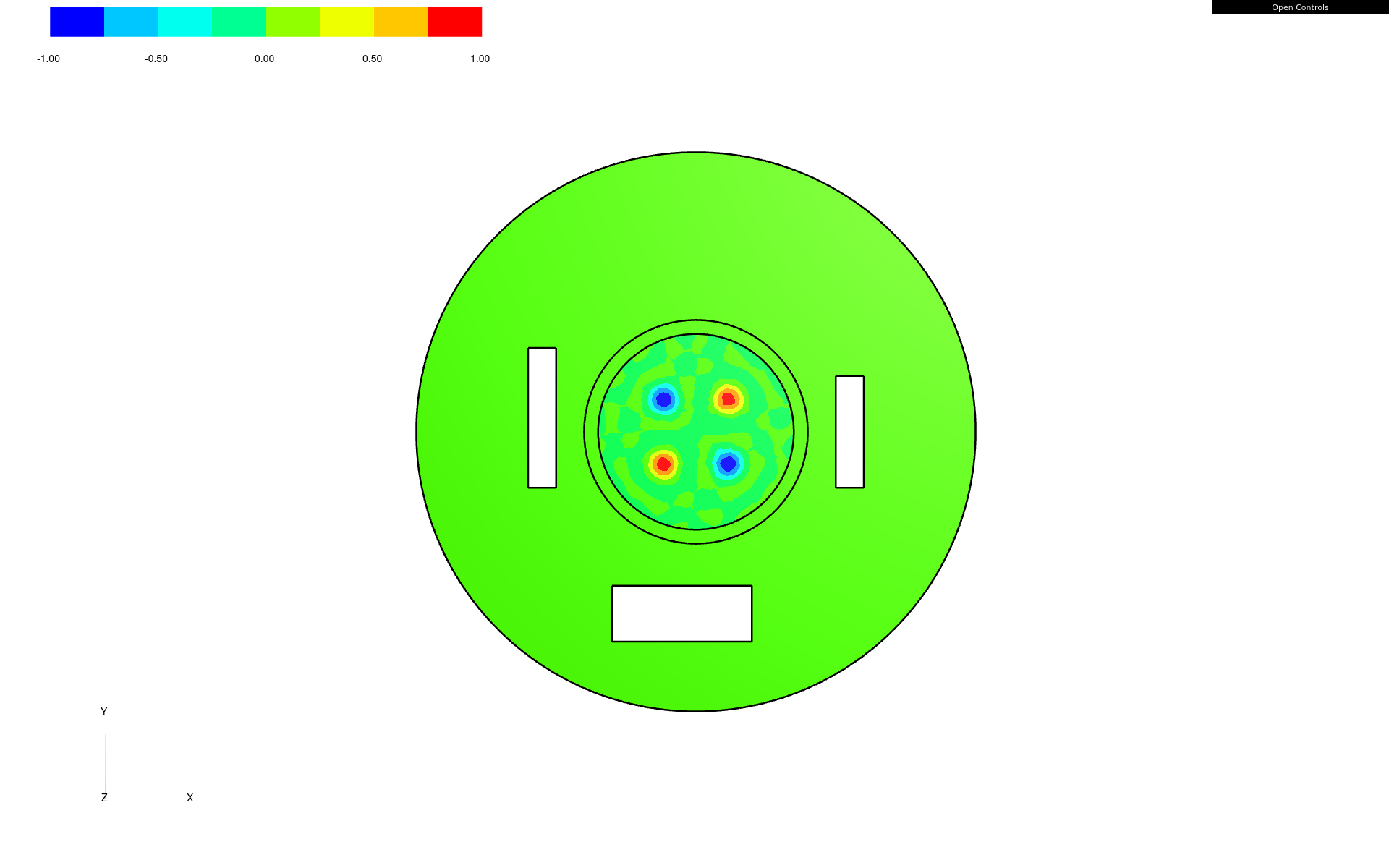} \\[-2ex]
\adjustimage{width=0.2\figsize,height=0.9\figsize,margin*={1.2cm 0.4cm 1.3cm -1cm}}{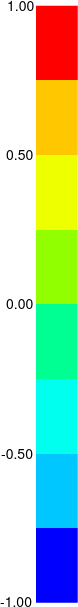}
\includegraphics[trim={19.8cm 5cm 19.8cm 6cm},clip,width=\figsize,keepaspectratio]{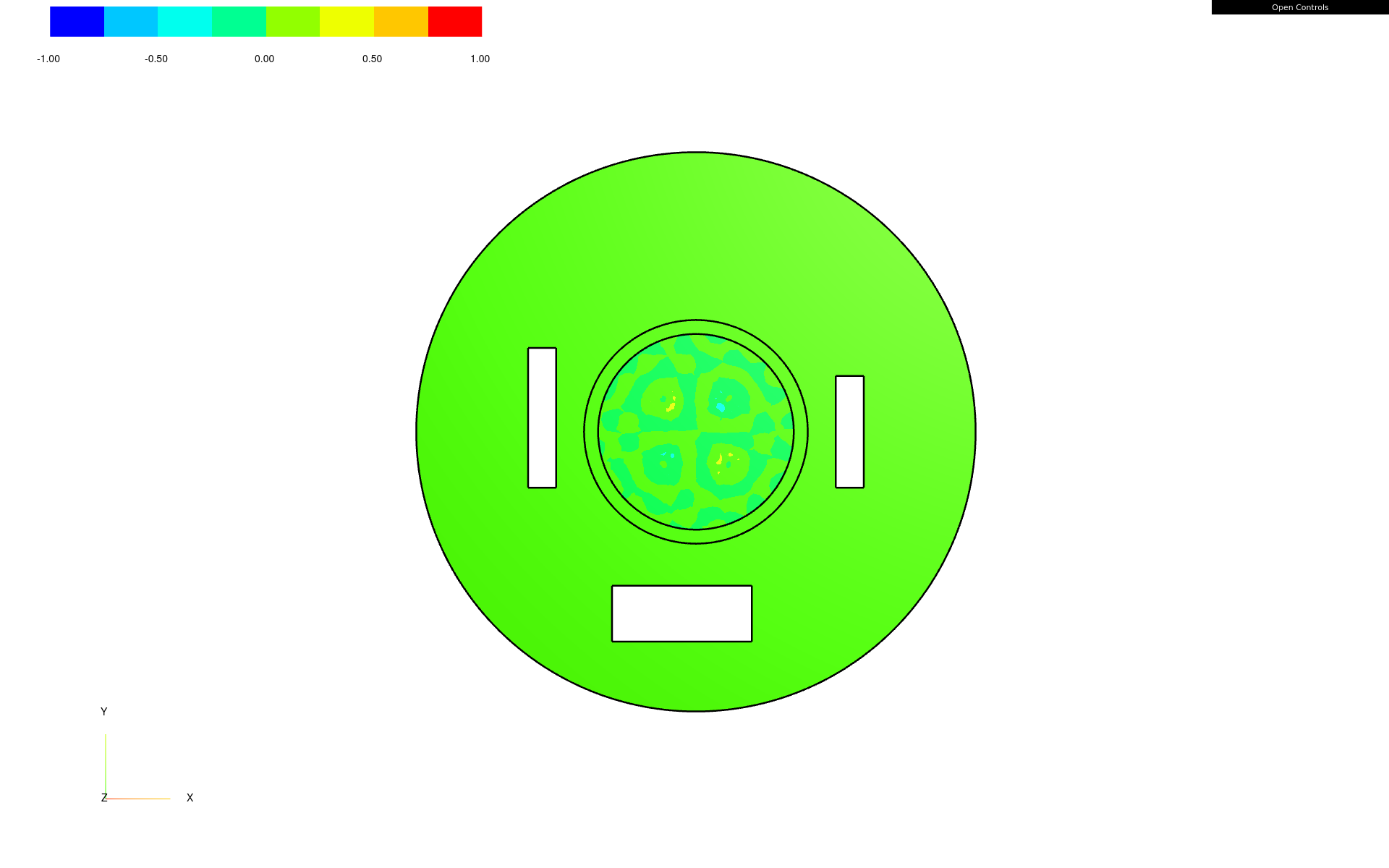}
\includegraphics[trim={19.8cm 5cm 19.8cm 6cm},clip,width=\figsize,keepaspectratio]{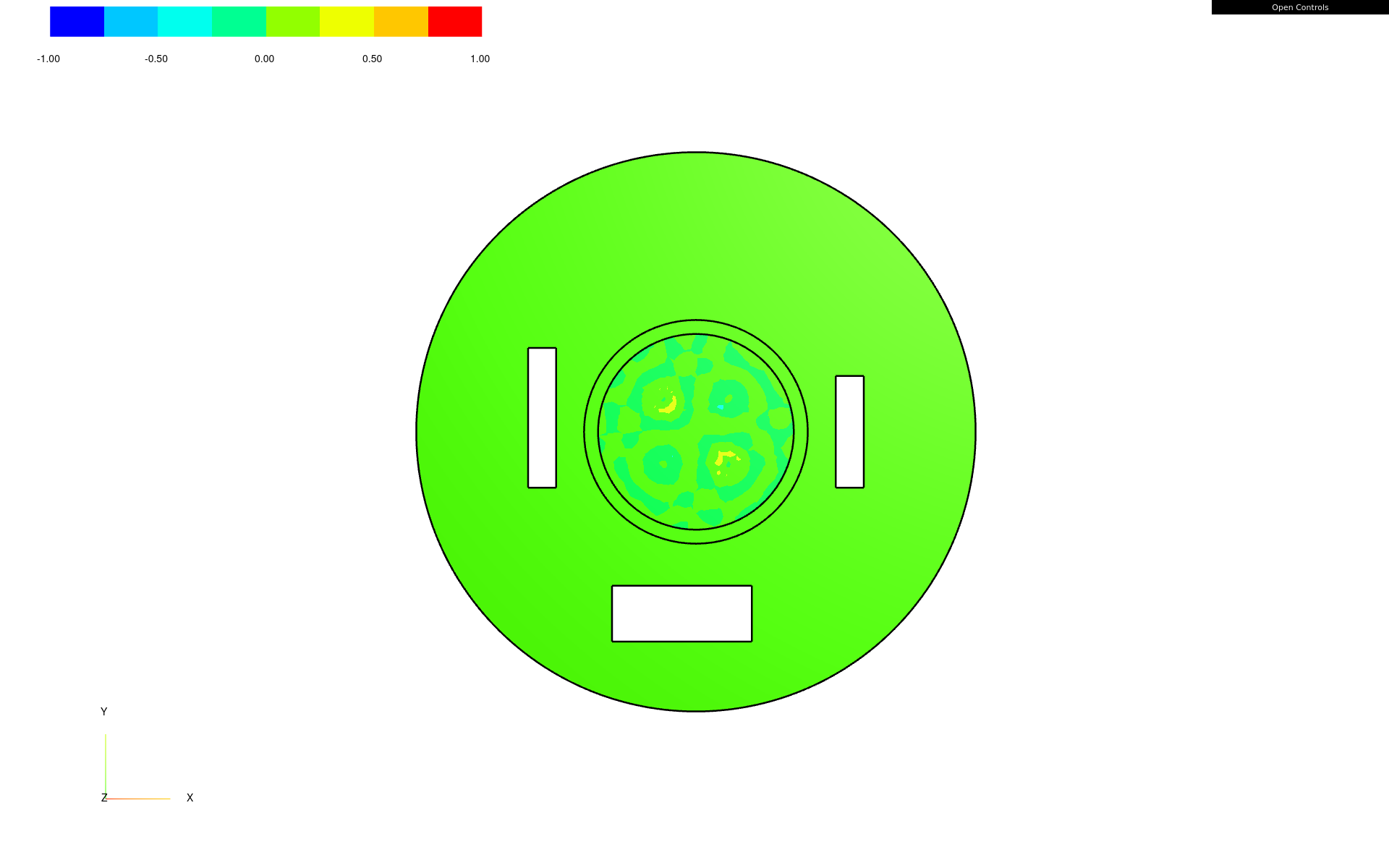}
\caption{Left: $f$. Middle column: Posterior mean, reconstruction error (global optimum-informed design, $m_0=12$). Right column: Posterior mean, reconstruction error (greedy design, $m_0=12$).}
\label{fig:recos}
\end{minipage}
\vspace{-0.5cm}
\end{figure}

On a 12th Gen Intel(R) Core(TM) i5-12500H (4.50 GHz) processor with 16 cores, Algorithm \ref{alg:greedy} returned a greedy sequence $(w^{m_0})_{m_0=1}^{36}$ over $100.0$ seconds. Meanwhile, foregoing the first five designs to allow the greedy algorithm to enter a sub-optimal domain, Algorithm \ref{alg:fgreedy} produced the global optimum-informed sequence $(w^{m_0}_*)_{m_0=6}^{36}$ over $701.4$ seconds, $629.1$ thereof computing the global optima $\ws$. While significantly slower than Algorithm \ref{alg:greedy}, the increased time yielded an average improvement of $4.2\%$ and a best-case improvement of $22.3\%$ in terms of the A-optimal objectives $\Jc$ compared to the greedy sequence; comparison can be seen in Figure \ref{fig:objectives}, showcasing that the global optimum-informed greedy sequence produced a better design for every value of $m_0$, and that it is significantly closer to the A-optimal objective value of the (non-binary) global optima $\ws$. Figure \ref{fig:recos} showcases posterior means $\mpst(w)$ for $m_0=12$ with the global optimum-informed resp.~greedy designs, given finite observations $F_wf$ of the source $f(x_1,x_2):=\sum_{i=0}^3(-1)^i\exp(-800\|(x_1-(-1)^{i+\delta_{i\geq 2}}r,x_2-(-1)^{\delta_{i\leq 1}}r)\|^2)$, $r:=0.35/3$. While both designs lead to good reconstructions, the greedy design had $L^2$ reconstruction error of approximately $4.69\cdot 10^{-2}$, while the global optimum-informed design had $L^2$ error of approximately $4.57\cdot 10^{-2}$, a roughly $2.7\%$ improvement.

\section{Conclusion and outlook}

We have demonstrated the power of the optimality criteria Theorem \ref{thm:optimality} for greedy designs, showcasing how it informs removal of underperforming sensors. Various improvements to Algorithm \ref{alg:fgreedy} can be made, as it does not take advantage of \emph{dominant} indices satisfying $\ws_k=1$, and the current greedy update step often re-introduces \emph{redundant} indices satisfying $\ws_k=0$. Efficiency-wise, performing the global-informed update step of Algorithm \ref{alg:fgreedy} less often might lead to significant speedup.

\begin{acknowledgement}
Work on this article was done with the support of the DFG's Grant 432680300 - SFB 1456 (C04). The author moreover expresses gratitude to Thorsten Hohage, University of Göttingen, and to Georg Stadler, Courant Institute of Mathematical Sciences.
\end{acknowledgement}




\end{document}